\documentclass[12pt]{article}
\usepackage{latexsym, amsmath, amssymb, amsfonts, amscd}
\usepackage{color}
\usepackage{graphics}
\usepackage[dvips]{epsfig}

\newtheorem{theorem}{Theorem}
\newtheorem{proposition}{Proposition}
\newtheorem{lemma}{Lemma}
\newtheorem{corollary}{Corollary}

\newtheorem{definition}{Definition}

\newcommand{\huaF}{{\cal F}}

\newcommand{\adL}{\mbox{\rm ad}_{\Lambda}}

\def\dbar{\overline\partial}

\newcommand{\CC}{\mathbb{C}}

\def\dbar{\overline\partial}

\def\partialbarc{{\overline{\partial}}_{\lie{c}}}
\def\partialbart{{\overline{\partial}}_{\lie{t}}}

\def\partialbar{\overline{\partial}}

\def\oomega{\overline\omega}
\def\Oomega{\overline\Omega}

\newcommand{\lie}[1]{\mathfrak{#1}}
\newcommand{\lieg}{\mathfrak{g}}
\newcommand{\liet}{\mathfrak{t}}
\newcommand{\liec}{\mathfrak{c}}

\newcommand{\bproof}{\noindent{\it Proof: }}
\newcommand{\eproof}{\hfill \qed \vspace{0.2in}}
\def\qed{\rule{2.3mm}{2.3mm}}

\begin{document}
\title{\bf Holomorphic Poisson Structures \\ and its Cohomology on Nilmanifolds}
\author{
Zhuo Chen\thanks{Address: Department of Mathematical Sciences, Tsinghua University, Beijing, P.R.C..
E-mail: zchen@math.tsinghua.edu.cn.} \ \
Anna Fino\thanks{ Address:  Dipartimento di Matematica G. Peano, Universit\`{a} degli Studi di Torino, Torino, Italy.
E-mail: annamaria.fino@unito.it} \ \
Yat-Sun  Poon\thanks{ Address:
    Department of Mathematics, University of California at Riverside,
    Riverside, CA 92521, U.S.A.. E-mail: ypoon@ucr.edu.} }

\date{}

\maketitle

\begin{abstract} The subject for investigation in this note is concerned with holomorphic Poisson
structures on nilmanifolds with abelian complex structures.
As a basic fact, we establish that on such manifolds, the Dolbeault cohomology with coefficients in
holomorphic polyvector fields is  isomorphic to the cohomology of invariant forms with
coefficients in invariant polyvector fields.

We then quickly identify the existence of invariant holomorphic Poisson structures. More important, the
spectral sequence of the Poisson bi-complex associated to such holomorphic Poisson structure degenerates at
$E_2$. We will also provide
examples of holomorphic Poisson structures on such manifolds so that the related spectral sequence
does not degenerate at $E_2$.

\end{abstract}


\section{Introduction}

The investigation of Poisson bracket from a complex perspective started a while ago \cite{Polish}.
Attention on this subject in the past ten years is largely due to its role in generalized complex
geometry   {\cite{Marco, Hitchin-Generalized CY, Hitchin-Instanton}}. It is now also known
that product of holomorphic Poisson structures together with symplectic structures forms the local
model of all generalized complex geometry \cite{Bailey}.

Therefore, one could take many different routes when investigating holomorphic Poisson structures.
One could study it as a complex geometric object and study its deformations as in
\cite{Hitchin-holomorphic Poisson}.
One could also investigate it in the context of extended deformation, or
generalized complex structures   {\cite{GPR, Poon}}.
In the former case, the deformation theory is dictated by the
differential Gerstenhaber structure on a cohomology ring with coefficients in the holomorphic polyvector
fields. In the latter case, it is on the part of the cohomology with total degree-2 as explained in
\cite{Poon}. As a key feature of generalized complex geometry is to encompass classical complex structures
with symplectic structures in a single geometric framework, one could also relate the cohomology theory
of a holomorphic Poisson structure as a generalized complex manifold to the cohomology theory on
symplectic geometry  {\cite{ACK, Tseng-Yau}}.

In this paper, we consider holomorphic Poisson structure as a geometric object in generalized complex
structure, and study its cohomology theory accordingly. It will set a stage for studying deformation theory.
In this perspective, it is known that the cohomology of a holomorphic Poisson structure could be computed by
a bi-complex   {\cite{Hong-Xu, Xu}}. The first level of the associated spectral sequence of this
bi-complex is the
Dolbeault cohomology with coefficients in holomorphic polyvector fields. It is known that
this spectral sequence often, but does not always degenerate at its second level \cite{CGP}.
It is therefore interesting to find how often this spectral sequence indeed degenerates
at its second level. On K\"ahlerian
manifolds, an affirmative answer for complex  {surfaces} is found, and other general observation is made
in \cite{CGP}.

In this note, we focus on non-K\"ahlerian manifolds. In particular, we focus on nilmanifolds due to its
rich history and role in generalized complex geometry \cite{CG}. Investigation on the cohomology theory
on nilmanifolds also has a very rich history, beginning with Nomizu's work on de Rham cohomology \cite{Nomizu}.
There has been a rich body of work on the Dolbeault cohomology of nilmanifolds with invariant complex
structures   {\cite{Console-Fino, CFGU}}, and work on Dolbeault cohomology on the same kind of manifolds
with coefficients in holomorphic tangent bundle  {\cite{Console, CFP, GMPP, MPPS-2-step,
Rolle}}. In
favorable situations, various authors proved that the cohomology is isomorphic to the cohomology
of invariant objects.

In this paper, after a review of holomorphic Poisson structures and their associated bi-complex structures and
a brief review of abelian complex structures on nilmanifolds,
we show that the Dolbeault cohomology of an abelian complex structure on a nilmanifold with
coefficients in holomorphic polyvector fields is  isomorphic to the cohomology of the
 corresponding invariant objects; see
Theorem \ref{quasi isomorphic s step}. It means that
the cohomology could be computed by a differential algebra over the field of
complex numbers. It enables an analysis of the  spectral sequence of the bi-complex associated to an
invariant holomorphic Poisson structure.

After we establish the existence of invariant holomorphic Poisson structures on nilmanifolds with abelian
complex structures in Section \ref{existence of Poisson}, we focus on proving Theorem \ref{existence}. This
theorem, which is also the key observation in this paper,
states that on any nilmanifold with abelian complex structures, there exists an invariant holomorphic
Poisson {structure} such that the spectral sequence of its associated bi-complex degenerates at its second level.
This result generalizes one of the observations in \cite{CGP} where the authors
could only work on 2-step nilmanifolds.

However, at the end of this note, we caution the readers
with an example that although such degeneracy occurs often, but it is not always true even in the context of
nilmanifolds with abelian complex structures.

\section{Holomorphic Poisson cohomology}

\

In this section, we review the basic background materials as seen in \cite{CGP} to set up the notations.

\
Let $M$ be a manifold with an integrable complex structure $J$.
Its complexified tangent bundle $TM_\CC$ splits into the direct sum of
bundle of $(1,0)$-vectors $TM^{1,0}$ and bundle of $(0,1)$-vectors $TM^{0,1}$.
Their $p$-th exterior products are respectively denoted by $TM^{p,0}$ and $TM^{0,p}$.
Denote their dual bundles by $TM^{*(p,0)}$ and $TM^{*(0,p)}$ respectively.

When $X,Y$ are vector fields, we denote their Lie bracket by $[X,Y]$. When $\omega$ is a 1-form,
we denote the Lie derivative of $\omega$ along $X$ by $[X,\omega]$. When $\rho$ is another 1-form,
we set $[\omega, \rho]=0$. With this \lq \lq bracket"  structure and the natural projection from
$ {\overline{L}}:= TM^{0,1}\oplus TM^{*(1,0)}$ to the summand $TM^{0,1}$, the bundle
${\overline{L}}$ is equipped with a complex Lie algebroid structure. Together with its conjugate bundle
$L= TM^{1,0}\oplus TM^{*(0,1)} $, they form a Lie bi-algebroid \cite{LWX}.
Then we get  the Lie algebroid differential $\dbar$ for the Lie algebroid $\overline{L}$ \cite{Mac}.
\begin{equation}
\dbar: C^\infty(M, TM^{1,0}\oplus TM^{*(0,1)} ) \to
C^\infty(M, \wedge^2( TM^{1,0}\oplus TM^{*(0,1)} )).
\end{equation}
It is extended to a differential of exterior algebras:
\begin{equation}
\dbar: C^\infty(M, \wedge^p(TM^{1,0}\oplus TM^{*(0,1)} )) \to
C^\infty(M, \wedge^{p+1}( TM^{1,0}\oplus TM^{*(0,1)} )).
\end{equation}
It is an exercise in Lie algebroid theory that the Lie algebroid differential
\[
\dbar: C^\infty(M, TM^{*(0,1)}) \to C^\infty(M, TM^{*(0,2)})
\]
is the $(0,2)$-component of the exterior differential, and
\[
\dbar:  C^\infty(M, TM^{1,0} )\to C^\infty(M, TM^{*(0,1)}\otimes TM^{1,0})
\]
is the Cauchy-Riemann operator \cite{Gau}.

On the space $C^\infty (M, \wedge^\bullet(TM^{1,0}\oplus TM^{*(0,1)} ))$, the Schouten bracket,
exterior product and the Lie algebroid differential $\dbar$ form a differential Gerstenhaber
algebra  {\cite{Mac,Poon}}.

Suppose that $\Lambda$ is a holomorphic Poisson structure, i.e. a smooth section of
$TM^{2,0}$ such that $[\Lambda, \Lambda]=0$ and $\dbar\Lambda=0$.
Denote the Schouten bracket of $\Lambda$ with elements in
$C^\infty (M, \wedge^\bullet(TM^{1,0}\oplus TM^{*(0,1)} ))$ by $\adL$, and the
 action of $\dbar+\adL$ on the same space by $\dbar_{\Lambda}$.
 Since $\Lambda$ is holomorphic Poisson,
 \begin{equation}
 \dbar_{\Lambda}: C^\infty (M, \wedge^k(TM^{1,0}\oplus TM^{*(0,1)} ))
 \to C^\infty (M, \wedge^{k+1}(TM^{1,0}\oplus TM^{*(0,1)} ))
 \end{equation}
 form an elliptic complex.
 From now on, for $n\geq 0$ denote
 \begin{equation}
 K^n=C^\infty (M, \wedge^n(TM^{1,0}\oplus TM^{*(0,1)} )).
 \end{equation}
 For $n< 0$, set $K^n=\{ 0\}$.
 \begin{definition}
For all $k\geq 0$, the $k$-th Poisson cohomology of the holomorphic Poisson structure $\Lambda$ is the space
 \begin{equation}
 H^k_{\Lambda}(M):=\frac{{\mbox{\rm kernel of  }}\ \dbar_\Lambda: K^k \to K^{k+1} }
 {{\mbox{\rm image of }}\ \dbar_\Lambda: K^{k-1}\to K^k}.
 \end{equation}
 \end{definition}
 Due to the nature of $\Lambda$, we have
 \begin{equation}
 \dbar\circ\dbar=0, \quad
 \dbar\circ \adL+\adL\circ \dbar=0,
 \quad
 \adL\circ\adL=0.
 \end{equation}
 In fact, the second identity is equivalent to $\Lambda$ being holomorphic, and
 the third is equivalent to $\Lambda$ being Poisson.
Define $A^{p,q}=C^\infty(M, TM^{p,0}\otimes TM^{*(0,q)})$, then
\begin{equation}
\adL: A^{p,q}\to A^{p+1, q},
\quad
\dbar: A^{p,q} \to A^{p, q+1};
\quad
\mbox{ and }
\quad
K^n=\oplus_{p+q=n}A^{p,q}.
\end{equation}
\begin{definition} Given a holomorphic Poisson structure $\Lambda$, the Poisson bi-complex is
the triple $\{ A^{p,q}, \adL, \dbar\}$.
\end{definition}
Then the cohomology $H^\bullet_\Lambda(M)$ theoretically could be computed by each one of the
two associated spectral sequences.
We choose a filtration given by $F^pK^n=\oplus_{p'+q=n, p'\geq p}A^{p',q}$. The
lowest differential is $\dbar: A^{p,q} \to A^{p, q+1}$. Therefore, the first level of the
spectral sequence is the Dolbeault cohomology
\begin{equation}
E_1^{p,q}=H^q(M, \Theta^p),
\end{equation}
where $\Theta^p$ is the sheaf of germs of the $p$-th exterior power of the holomorphic tangent bundle
on the complex manifold $M$. It follows that the next differential is
 \begin{equation}
 d_1^{p,q}=\adL: H^q(M, \Theta^p)\to H^q(M, \Theta^{p+1}).
 \end{equation}
The second level of the Poisson spectral sequence is given by
\begin{equation}\label{Eqt:E2pq}
E_2^{p,q}=
\frac{{\mbox{\rm kernel of  }}\ \adL:H^q(M, \Theta^p)\to H^q(M, \Theta^{p+1})}
{{\mbox{\rm image of  }}\ \adL:H^q(M, \Theta^{p-1})\to H^q(M, \Theta^{p})}.
\end{equation}
We are interested in computing
\begin{equation}
d_2^{p,q}: E_2^{p,q}\to E_2^{p+2,q-1}.
\end{equation}

\section{Nilmanifolds with abelian complex structures}\label{Zhuo}

A compact manifold $M$ is called a nilmanifold if there exists a simply-connected nilpotent Lie group $G$
and a lattice subgroup $\Gamma$ such that $M$ is diffeomorphic to $G/\Gamma$.
We denote the Lie algebra of
$G$ by $\lieg$. The step of the
nilmanifold is the nilpotence of the Lie algebra $\lieg$. For a $(s+1)$-step nilpotent Lie algebra $\lieg$, there is
a filtration by the descending central series,
\[
\left\{  0\right\}  =\lie{g}^{s+1}=
[\mathfrak{g}^{s}, \lie{g}]\subset
\cdots
\subset\mathfrak{g}^{k+1}
=
[  \mathfrak{g}^{k}, \mathfrak{g}]
\subset\mathfrak{g}^{k}\subset\cdots\subset
\mathfrak{g}^{1}
=[  \mathfrak{g}, \lie{g}]  \subset
\lie{g}^{0}=\lie{g}.
\]

A left-invariant complex
structure $J$ on $G$ is said to be abelian if on the Lie algebra $\lieg$, it satisfies the condition
$[JA,JB]=[A,B],$
for all $A$ and $B$ in the Lie algebra $\lieg$.
If one complexifies the algebra $\lieg$ and denotes the $+i$ and $-i$ eigen-spaces of $J$ respectively
by $\lieg^{1,0}$ and $\lieg^{0,1}$, then the invariant complex structure $J$ being abelian is equivalent
to the complex algebra $\lieg^{1,0}$ being abelian.

Denote the $k$-th exterior products $\wedge^k\lieg^{1,0}$ and $\wedge^k{\lieg}^{*(0,1)}$ respectively by
$\lieg^{k,0}$ and ${\lieg}^{*(0,k)}$. There is a natural inclusion map.
\[
\iota:{\lie g}^{\ell,0}\otimes {\lie g}^{*(0,m)}
\rightarrow C^\infty(M, \wedge^\ell TM^{1,0} \otimes \wedge^m TM^{*(0,1)} ).
\]

Now we wish to prove the following:

 \begin{theorem}\label{quasi isomorphic s step} On a nilmanifold $M$ with an invariant
abelian complex structure, the inclusion ${\lie g}^{\ell,0}\otimes {\lie g}^{*(0,m)}$ in
$C^\infty(M, \wedge^\ell TM^{1,0} \otimes \wedge^m TM^{*(0,1)} )$
 induces an isomorphism of cohomology. In other words,
$H^m({\lieg}^{\ell,0})\cong H^m(M, \Theta^\ell)$.
\end{theorem}
This theorem generalizes an observation in \cite{CGP} for 2-step nilmanifolds.

\section{The proof of Theorem \ref{quasi isomorphic s step}}

\

We adopt an inductive approach regarding the number of steps because the
Theorem \ref{quasi isomorphic s step} is known to be true for 2-step  nilmanifolds
 \cite{CGP}.

Let   $\lieg$ be $(s+1)$-step nilpotent.  Assume that the theorem holds true for $s$-step nilmanifolds where $s\geq 1$.

 Let $\lie{c}$ be the center of the Lie algebra $\lie{g}$. Since
$[JA, B]=-[A, JB]$, the center $\lie{c}$ is $J$-invariant.
Let $\lie{t}=\lie{g}/\lie{c}$.
It is obvious that $\lie{t}$ is $s$-step nilpotent
and it has an induced abelian complex structure as well.

Let $C$ be the center of $G$ and $\psi: G\rightarrow G/C$ the quotient map.
Since $G$ is $(s+1)$-step nilpotent, $G/C$ is $s$-step nilpotent.
Consider $M=G/\Gamma$ and $N=\psi(G)/\psi(\Gamma)$.
We have a holomorphic fibration $\Psi: M\rightarrow N$ whose fiber is isomorphic to $F=C/(C\cap \Gamma)$.
Note also that $N$ is a $s$-step nilmanifold with an abelian complex structure.

We have   {the} vector space decompositions
\[
{\lie{g}}_{\mathbb{C}}={\lie{g}}^{1,0}\oplus {\lie{g}}^{0,1};
\quad
{\lie{g}}^{1,0}={\lie{c}}^{1,0}\oplus {\lie{t}}^{1,0};
\quad
{\lie{g}}^{0,1}={\lie{c}}^{0,1}\oplus {\lie{t}}^{0,1}.
\]
Here both ${\lie{g}}^{1,0}$ and ${\lie{g}}^{0,1}$ are abelian sub-algebras of
$\lieg_{\CC}$. The only non-trivial Lie brackets are of the form:
\[
[\lie{t}^{1,0},\lie{t}^{0,1}]\subset \lieg_{\CC}= \lie{c}^{1,0}\oplus \liec^{0,1}
\oplus  \lie{t}^{1,0}\oplus \lie{t}^{0,1}.
\]
For the $d=\partial + \dbar $ operator, we have the following lemma, which can be verified directly.
\begin{lemma} We have $\dbar \lieg^{*(0,1)}=0$,
$\dbar \lie{g}^{*(1,0)}\subset \lie{t}^{*(1,1)}$, $\dbar  \lie{c}^{1,0}=0$ and
\begin{equation}\label{general dbar}
\dbar  \lie{t}^{1,0}\subset
(\lie{t}^{*(0,1)}\otimes \lie{c}^{1,0}) ~\oplus~ (\lie{t}^{*(0,1)}\otimes \liet^{1,0}).
\end{equation}
\end{lemma}

To compute the cohomology $H^m(\lieg^{\ell,0})$, one uses
the $\dbar$-operator:
\[
\dbar :~  \lie{g}^{*(0,m)}\otimes \lie{g}^{\ell,0}\rightarrow \lie{g}^{*(0,m+1)}\otimes \lie{g}^{\ell,0}.
\]
As $\lie{g}^{1,0}=\lie{c}^{1,0}\oplus \lie{t}^{1,0}$, we have
\[
\lie{g}^{\ell,0}=\oplus_{a+b=\ell} \lie{c}^{a,0}\otimes \lie{t}^{b,0};
\qquad \lie{g}^{*(0,m)}=\oplus_{i+j=m} \lie{c}^{*(0,i)}\otimes \lie{t}^{*(0,j)}.
\]
According to the {decomposition}  in Equation (\ref{general dbar}),
one may split the $\partialbar $ operator into two parts,
\[
\dbar =\dbar_{\lie{c}}+\dbar_{\lie{t}},
\]
depending on whether we choose the $\liec^{1,0}$ component or $\liet^{1,0}$ in the range of the
operator $\dbar$.
Here
\begin{equation}\label{Eqt:dc}
\partialbarc:~\liec^{*(0,i)}\otimes \liet^{*(0,j)}\otimes \liec^{a,0}\otimes \liet^{b,0}
\rightarrow \liec^{*(0,i)}\otimes \liet^{*(0,j+1)}\otimes \liec^{a+1,0}\otimes { \liet^{b-1,0}}
\end{equation}
and
\begin{equation}\label{Eqt:dt}
\partialbart:~\liec^{*(0,i)}\otimes \liet^{*(0,j)}\otimes \liec^{a,0}\otimes \liet^{b,0}
\rightarrow \liec^{*(0,i)}\otimes \liet^{*(0,j+1)}\otimes \liec^{a,0}\otimes  {\liet^{b,0}}.
\end{equation}

We now fix the number $\ell\geq 0$. If for $0\leq p \leq \ell$, we put
\[
D^{p,q}= \oplus_{i+j=p+q} ~~\liec^{*(0,i)}\otimes \liet^{*(0,j)}
\otimes \liec^{p,0}\otimes \liet^{\ell-p,0}
=\oplus_{0\leq i\leq p+q} ~~\liec^{*(0,i)}\otimes \liet^{*(0,p+q-i)}
\otimes \liec^{p,0}\otimes \liet^{\ell-p,0},
\]
then
\[
\partialbarc: D^{p,q}\rightarrow D^{p+1,q},
\quad
\mbox{ and }
\quad
\partialbart: D^{p,q}\rightarrow D^{p,q+1}.
\]
And hence
\[
\partialbarc\circ\partialbarc: D^{p,q}\rightarrow D^{p+2,q},
\quad
\partialbarc\circ\partialbart
+\partialbart\circ\partialbarc:
D^{p,q}\rightarrow D^{p+1,q+1},
\quad
\partialbart\circ\partialbart:
D^{p,q}\rightarrow D^{p,q+2}.
\]
Since $\dbar=\partialbarc+\partialbart$, and
$\dbar\circ\dbar=0$, we have
\[
\partialbarc\circ\partialbarc=0,
\quad
\partialbarc\circ\partialbart
+\partialbart\circ\partialbarc=0,
\quad
\partialbart\circ\partialbart=0.
\]
It means that both $\partialbarc$ and $\partialbart$ are co-boundary operators, and the data
$(D^{p,q},\partialbarc,\partialbart
)$ form a bi-complex.
Its total complex is
\[
\partialbar=\partialbarc+\partialbart: \oplus_{p+q=m}D^{p,q}\rightarrow \oplus_{p+q=m+1}D^{p,q},
\]
or, exactly that of
\[
\lieg^{*(0,m)}\otimes \lieg^{\ell,0}\rightarrow \lieg^{*(0,m+1)}\otimes
\lieg^{\ell,0}.
\]

In summary, we have
\begin{lemma} For each $\ell\geq 0$, the cohomology $H^m(\lieg^{\ell,0})$ can be computed
as the total cohomology of the bi-complex $(D^{p,q},\partialbarc,\partialbart)$.
\end{lemma}
\begin{lemma}\label{Lemma:intermiedietfacts} For each $\ell \geq 0$, we have
\[H^m(M,\wedge^\ell\Psi^* \Theta_N )=\oplus_{p+q=m} \liec^{*(0,q)}\otimes H^p(\liet^{\ell,0}) .\]
\end{lemma}
\bproof The proof is essentially the same as that of Lemma 5 in \cite{MPPS-2-step}.
Here we give a sketch.

To compute $H^m(M,\wedge^\ell\Psi^* \Theta_N )$, we use the standard
Leray spectral sequence for a fibration. We first note the following fact (see Lemma 3 in \cite{MPPS-2-step}):
\[
R^q\Psi_*(\wedge^\ell \Psi^*\Theta_N)=\liec^{*(0,q)}\otimes \wedge^\ell\Theta_N.
\]
Therefore, the second level of the Leray spectral sequence is given by
\[
E^{p,q}_2 = H^p(N,R^q\Psi_*(\wedge^\ell\Psi^*\Theta_N))=H^p(N,\liec^{*(0,q)}\otimes \wedge^\ell\Theta_N)
=\liec^{*(0,q)}\otimes H^p(N, \wedge^\ell\Theta_N).
\]
By induction assumptions: $N$ is a $s$-step nilmanifold and
$H^p(N,  \wedge^\ell\Theta_N)=H^p(\liet^{\ell,0})$. Therefore,
\begin{equation}
\liec^{*(0,q)}\otimes H^p(N, \wedge^\ell\Theta_N)=\liec^{*(0,q)}\otimes H^p(\liet^{\ell,0}).
\end{equation}

Thus $d_2: E^{p,q}_2\rightarrow E^{p+2,q-1}_2$ is a map
\[
\liec^{*(0,q)}\otimes   H^p(\liet^{\ell,0})\rightarrow
\liec^{*(0,q-1)}\otimes   H^{p+2}(\liet^{\ell,0}).
\]
However, any element  in $\liec^{*(0,q)}\otimes   H^p(\liet^{\ell,0})$ could be represented
in the following form:
\[
\sum_{a, b, c}{\overline{\rho}}_a\wedge \Oomega_b\wedge W_c\,,
\]
where ${\overline{\rho}}_a\in \liec^{*(0,q)}$,
$\Oomega_b\in \liet^{*(0,p)}$, $W_c\in \liet^{\ell, 0}$, and
 $\dbar_N(\sum_{b, c} \Oomega_b\wedge W_c)=0$.
 Note that for forms, $\dbar_N\Oomega_b=\dbar_M\Oomega_b=0$ because the complex structures on both the
 manifold $M$ and its quotient $N$ are abelian.
 In addition, since the fibers of the projection $\Psi$ are global holomorphic vector fields generated by
 $\liec^{1,0}$, and they are in the center of $\lie{g}_\CC$,
 \[
 \dbar_N W_c=\dbar_M W_c.
 \]
 Therefore, when $\dbar_N(\sum_{b, c} \Oomega_b\wedge W_c)=0$, then
 \[
 \dbar \left(\sum_{a, b, c}{\overline{\rho}}_a\wedge \Oomega_b\wedge W_c \right)=0.
 \]
 It follows that $d_2\equiv 0$, and
\begin{align*}
H^m(M,\wedge^{\ell}\Psi^* \Theta_N )&=\oplus_{p+q=m}  E^{p,q}_2
=\oplus_{p+q=m}  ~\liec^{*(0,q)}\otimes  H^p(\liet^{\ell,0}).
\end{align*}
\eproof

Now we are back to the proof of Theorem 1.
We take   $E=\Theta_M$,
$Z=\mathcal{O}_M\otimes \liec^{1,0}  $, $Q=\Psi^*\Theta_N$.
Then we have the following exact sequence of holomorphic vector bundles over $M$:
\[
0\rightarrow \mathcal{O}_M\otimes \liec^{1,0} \rightarrow
\Theta_M \stackrel{\rho}{\longrightarrow} \Psi^*\Theta_N \rightarrow 0,
\]
{i.e.},
\[
0\rightarrow Z\rightarrow E\stackrel{\rho}{\longrightarrow} Q\rightarrow 0.
\]

By Lemma 5 in \cite{CGP}, we have a filtration of $\wedge^\ell E=\wedge^\ell \Theta_M$:
\[
\mathcal{O}_M\otimes \liec^{\ell,0}=\wedge^{\ell} Z =E^{(\ell)}\subset E^{(\ell-1)}
\subset \cdots \subset E^{(1)}\subset E^{(0)}=\wedge^\ell E=\wedge^\ell \Theta_M.
\]
Moreover, the associated graded spaces are:
\[
G^0= E^{(0)} / E^{(1)}\cong \wedge^\ell Q=\Psi^*(\wedge^{\ell} \Theta_N);
\]
\[
G^1=E^{(1)} / E^{(2)}\cong  \wedge^{\ell-1} Q\otimes Z=\Psi^*(\wedge^{\ell-1} \Theta_N)\otimes \liec^{1,0};
\]
\[\cdots\cdots\cdots
\]
\[
G^r=  E^{(r)} / E^{(r+1)} =  \wedge^{\ell-r} Q\otimes \wedge^r Z
= \Psi^*(\wedge^{\ell-r} \Theta_N)\otimes \liec^{r,0};
\]
\[
\cdots\cdots\cdots
\]
\[G^{\ell-1}=E^{(\ell-1)} / E^{(\ell)}\cong  Q\otimes \wedge^{\ell-1}Z
= \Psi^*(  \Theta_N)\otimes \liec^{\ell-1,0}.
\]

Accordingly, we have a filtration of the co-chain complex
$C^\bullet = TM^{*(0,\bullet)}\otimes \wedge^\ell \Theta_M$:
\[
  TM^{*(0,\bullet)}\otimes \liec^{\ell,0}
=\huaF^{(\ell)} C^{\bullet}\subset \huaF^{(\ell-1)} C^{\bullet} \subset \cdots
\subset\huaF^{(1)} C^{\bullet}\subset \huaF^{(0)} C^{\bullet}= C^\bullet,
\]
where $\huaF^{(r)}C^\bullet=TM^{*(0,\bullet)}\otimes E^{(r)}$.

Thus there associates a spectral sequence: $E_r^{p,q}$, which starts with
$$
E^{p,q}_0= \huaF^{(p)}C^{p+q}/ \huaF^{(p+1)}C^{p+q}=
TM^{*(0,p+q)}\otimes G^p
\cong TM^{*(0, p+q)}\otimes \Psi^*(\wedge^{\ell-p} \Theta_N)\otimes \liec^{p,0}
$$
and $d_0=\partialbar $. It follows from Lemma
\ref{Lemma:intermiedietfacts} that we have
\[
 E^{p,q}_1=
H^{p+q}(G^p)=H^{p+q}( \Psi^*(\wedge^{\ell-p} \Theta_N)\otimes
\liec^{p,0})=\oplus_{i+j=p+q} \liec^{*(0,i)}\otimes H^j(\liet^{\ell-p,0}) \otimes\liec^{p,0}.
\]
It can be easily seen that the right hand side is in fact the cohomology of
\[
\partialbart:~ \oplus_{i+j=p+q} \liec^{*(0,i)}\otimes \liet^{*(0,j)}\otimes
\liet^{\ell-p,0} \otimes\liec^{p,0}
\rightarrow \oplus_{i+j=p+q} \liec^{*(0,i)}\otimes
\liet^{*(0,j+1)}\otimes \liet^{\ell-p,0} \otimes\liec^{p,0}.
\]
Thus,
\[
E^{p,q}_1=\frac{{\mbox{\rm kernel of }}\ \partialbart: ~D^{p,q}\rightarrow D^{p,q+1} }
{\mbox{\rm image of }\ \partialbart:~D^{p,q-1}\rightarrow D^{p,q}}.
\]

We then find the associated $d_1: E^{p,q}_1\rightarrow E^{p+1,q}_1$.
In fact, it is essentially $\partialbar$. If $X\in D^{p,q}$ represents
an element in $E^{p,q}_1$, then $d_1[X]$ is represented by $\partialbar X$.
Note $\partialbar=\partialbart+\partialbarc$
and $\partialbart X=0$. So $d_1[X]$ is actually represented by $\partialbarc X$.

Using this description, we now explain
\[
E^{p,q}_2
=\frac{{\mbox{\rm kernel of }}\ d_1: ~E^{p,q}_1\rightarrow E^{p+1,q}_1}
{{\mbox{\rm image of }}\ d_1:~E^{p-1,q}_1\rightarrow E^{p,q}_1} .
\]

An element in  $E^{p,q}_2$ can be represented by some $X\in D^{p,q}$
satisfying the following condition:~$\partialbart X=0$, and
$\exists$ $Y\in D^{p+1,q-1}$  such that $\partialbarc X+\partialbart Y=0$.

Moreover, such an $X$ represents the zero element in $E^{p,q}_2$ if
there exist  $Z\in D^{p-1,q}$ and  $Z'\in D^{p,q-1}$ such that $\partialbart Z=0$
and $\partialbarc Z+\partialbart Z'=X$.

The co-boundary at level $2$, $d_2: E^{p,q}_2\rightarrow E^{p+2,q-1}_2$,
which is again essentially $\partialbar$, now becomes the map sending $[X]$
to $[\partialbarc Y]$. In fact, this can be easily seen from the following calculation:
\begin{eqnarray*}
 \partialbar X  &=& \partialbart X+\partialbarc X=\partialbarc X\\
&=& - \partialbart Y \equiv \partialbarc Y \mod {\mbox{Im}}(\partialbar).
\end{eqnarray*}

Repeat this process. We find that $E^{p,q}_r=Z^{p,q}_r/B^{p,q}_r$, where
 $
Z^{p,q}_r
 $ consists of elements $X^{p,q}\in D^{p,q}$ which is subject to the following conditions:

 \begin{enumerate}
 \item $\partialbart X^{p,q}=0$;
 \item $\exists X^{p+1,q-1}\in  D^{p+1,q-1}$, $X^{p+2,q-2}\in  D^{p+2,q-2}$,
 $\cdots$, $X^{p+r-1,q-r+1}\in  D^{p+r-1,q-r+1}$ such that
     \begin{align*}
     &\partialbarc X^{p,q}+\partialbart X^{p+1,q-1} = 0,\\
     &\cdots  \cdots,\\
     &\partialbarc X^{p+r-2,q-r+2}+\partialbart X^{p+r-1,q-r+1} = 0.
     \end{align*}
 \end{enumerate}

The denominator $
B^{p,q}_r
 $ consists of elements $W^{p,q}\in D^{p,q}$  {satisfying}  the following condition:
 $\exists$ $W^{p,q-1}\in D^{p,q-1}$, $\cdots$, $W^{p-r+1,q+r-2}\in D^{p-r+1,q+r-2}$, such that
 \begin{align*}
 &\partialbart W^{p,q-1}+ \partialbarc W^{p-1,q}=W^{p,q},　\\　
  &\partialbart W^{p-1,q}+ \partialbarc W^{p-2,q+1}=0,　\\　
  &\cdots \cdots,\\
  &\partialbart W^{p-r+2,q+r-3}+ \partialbarc W^{p-r+1,q+r-2}=0,　\\
  &\partialbart W^{p-r+1,q+r-2}=0.
 \end{align*}

With this construction of $E^{p,q}_r=Z^{p,q}_r/B^{p,q}_r$,   the co-boundary map
$d_r: E^{p,q}\rightarrow E^{p+r,q-r+1}$ should be given by $d_r[X^{p,q}]=[\partialbarc X^{p+r-1,q-r+1}]$.
The reason is similar to the preceding $r=2$ case.

From these observations, we find that the spectral sequence $E^{p,q}_r$ derived from the filtration
of $\wedge^\ell E=\wedge^\ell \Theta_M$ exactly matches with
the spectral sequence of the bi-complex $(D^{p,q},\partialbarc,\partialbart)$.
Hence they converge to the same cohomology, i.e.,
\[
H^m(M,  \wedge^\ell\Theta_M)\cong H^m_{\mathrm{total}}(D^{p,q})\cong H^m(\lieg^{\ell,0}).
\]
\eproof

According to Theorem \ref{quasi isomorphic s step} and Equation \eqref{Eqt:E2pq},
\begin{equation}
E_2^{p,q}=
\frac{{\mbox{\rm kernel of }} \adL:H^q(\lieg^{p, 0})\to H^q(\lieg^{p+1})}
{{\mbox{\rm image of }} \adL:H^q(\lieg^{p-1})\to H^q(\lieg^{p})}.
\end{equation}

In this context, it is clear   {that} if the $A$ and $B$ are complex linearly independent elements
in $\liec^{1,0}$ and $\Lambda=A\wedge B$, then $\Lambda$ is a holomorphic Poisson structure
such that $\adL=0$. In such case, $\dbar_\Lambda=\dbar$ and
\begin{equation}
H^k_{\Lambda}=\oplus_{p+q=k}E^{p,q}_2=\oplus_{p+q=k}H^q(\lieg^{p,0}).
\end{equation}

Therefore, the only non-trivial part of Theorem \ref{existence} below is when the dimension of
$\liec^{1,0}$ is equal to one.

\begin{theorem}\label{existence}
 On any nilmanifold with an abelian complex structure, there exists a
 non-trivial holomorphic Poisson
structure $\Lambda$ such that its associated spectral sequence {degenerates} on its second level. In particular,
\begin{eqnarray}
H^k_{\Lambda} &=&\oplus_{p+q=k}E_2^{p,q}\\ \nonumber
&=&\oplus_{p+q=k}
\frac{\mbox{\rm kernel of } \adL:  H^q({\lieg}^{p,0})
 \to  H^q({\lieg}^{p+1,0})}{\mbox{\rm image of }
 \adL:  H^q({\lieg}^{p-1,0}) \to
  H^q({\lieg}^{p,0})}.
  \end{eqnarray}
\end{theorem}

\section{Existence of holomorphic Poisson structures}\label{existence of Poisson}

We will continue to work with abelian complex structures. Recall that by definition,
$\lieg^0=\lieg$, and inductively, $\lieg^{k+1}=[\lieg^k, \lieg]$.

Define $\mathfrak{g}_{J}^{k}=\mathfrak{g}^{k}+J\mathfrak{g}^{k}.$ When the complex structure
is abelian, it is clear from various definitions that each
$\mathfrak{g}_{J}^{k}$ is a $J$-invariant ideal of $\mathfrak{g}$, and
  we have a  filtration of subalgebras:
\[
\left\{  0\right\}  =\mathfrak{g}_{J}^{s+1}\subset\mathfrak{g}_{J}^{s}
\subset\mathfrak{g}_J^{s-1}
\subseteq\cdots\subseteq\mathfrak{g}_{J}^{k+1}\subseteq\mathfrak{g}%
_{J}^{k}\subseteq\cdots\subseteq\mathfrak{g}_{J}^{1} {\subset} \mathfrak{g}%
_{J}^{0}=\mathfrak{g}.
\]
Note that   by \cite{Salamon}  the last inclusion  $\mathfrak{g}_{J}^{1} \subset  \mathfrak{g}$ is always strict.    Moreover,  since the center $\liec$ is $J$-invariant and it contains $\lieg^s$,  we have
\begin{equation}\label{center}
\lie{g}^s_J\subseteq \liec.
\end{equation}
It follows that the   inclusion $\mathfrak{g}_{J}^{s}
\subset\mathfrak{g}_J^{s-1}$ is also strict.

We complexify this filtration:
\begin{equation}\label{filtration}
\left\{  0\right\}  =\mathfrak{g}_{J,\CC}^{s+1}\hookrightarrow\mathfrak{g}_{J,\CC}
^{s}\hookrightarrow\cdots\hookrightarrow\mathfrak{g}_{J,\CC}^{k+1}\hookrightarrow
\mathfrak{g}_{J,\CC}^{k}\hookrightarrow\cdots\hookrightarrow\mathfrak{g}_{J,\CC}^{1}\hookrightarrow
\mathfrak{g}_{J,\CC}^{0}=\mathfrak{g}_{\CC}.
\end{equation}
There exists a type decomposition for each $k$.
\[
\mathfrak{g}_{J,\CC}^{k}=\mathfrak{g}_{J}^{k,(1,0)}\oplus\mathfrak{g}_{J}^{k,(0,1)}.
\]
So, the filtration (\ref{filtration}) splits into two.
One is for type $(1,0)$-vectors.
\[
\left\{  0\right\}  \hookrightarrow\mathfrak{g}_{J}^{s,(1,0)}
\hookrightarrow\cdots\hookrightarrow\mathfrak{g}_{J}^{k+1,(1,0)}
\hookrightarrow\mathfrak{g}_{J}^{k,(1,0)}\hookrightarrow\cdots\hookrightarrow\mathfrak{g}_{J}^{1,(1,0)}
\hookrightarrow\mathfrak{g}_{J}^{0,(1,0)}=\mathfrak{g}^{(1,0)};
\]
Another is for type $(0,1)$-vectors.
\[
\left\{  0\right\}  \hookrightarrow\mathfrak{g}_{J}^{s,(0,1)}
\hookrightarrow\cdots\hookrightarrow\mathfrak{g}_{J}^{k+1,(0,1)}
\hookrightarrow\mathfrak{g}_{J}^{k,(0,1)}\hookrightarrow\cdots\hookrightarrow\mathfrak{g}_{J}^{1,(0,1)}
\hookrightarrow\mathfrak{g}_{J}^{0,(0,1)}=\mathfrak{g}^{(0,1)}.
\]

\begin{lemma}\label{lemma on brackets} Suppose that the complex structure $J$ is abelian, then
\begin{itemize}
\item
$\left[  \mathfrak{g}^{(1,0)},\mathfrak{g}^{(1,0)}\right] =0,$ and
$\left[  \mathfrak{g}^{(0,1)},\mathfrak{g}^{(0,1)}\right] =0.$
\item
$\left[  \mathfrak{g}_{J}^{k,(1,0)},\mathfrak{g}_{J}^{\ell,(0,1)}\right]
\subseteq\mathfrak{g}_{J, \CC}^{1+\max{\{k,\ell\}}}.$
\item In particular, when $\ell=0$,
$\left[  \mathfrak{g}_{J}^{k,(1,0)},\mathfrak{g}^{(0,1)}\right]
\subseteq\mathfrak{g}_{J, \CC}^{k+1}$
\end{itemize}
\end{lemma}
\bproof The first point is due to the complex structure being abelian. To prove the second point,
assume that $k\geq \ell$, let  $X$ be in $\lie{g}^k$ and $Y$ be in
 $\lie{g}^\ell$.
\begin{eqnarray*}
[X-iJX, Y+iJY]&=&[X,Y]+[JX,JY]-i([JX,Y]-[X,JY])\\
&=&2[X,Y]+2i[X,JY].
\end{eqnarray*}
As $k\geq \ell$, by definition $[X,Y]\in \lieg^{k+1}$ and $[X,JY]\in \lieg^{k+1}$. Therefore,
$[X-iJX, Y+iJY]$ is contained in $\lie{g}^{k+1}_\CC$.

In general, if $X_1$ and $X_2$ are in $\lie{g}^k$, then
\[
 X_1+JX_2-iJ(X_1+JX_2)=(X_1-iJX_1)+i(X_2-iJX_2).
 \]
 By complex linearity, the proof of the second observation is completed.
\eproof

Make the following notation for the quotient space
\[
\mathfrak{t}^{k+1,(1,0)}=\lie{g}_{J}^{k,(1,0)}/
\lie{g}_{J}^{k+1,(1,0)}.
\]


Choose a vector space isomorphism so that the short exact sequence of Lie algebras
\[
0\to \lie{g}_{J}^{k+1,(1,0)}\to \lie{g}_{J}^{k,(1,0)} \to \lie{g}_{J}^{k,(1,0)}/
\lie{g}_{J}^{k+1,(1,0)}\to 0
\]
is turned into a direct sum of vector spaces.
\[
\mathfrak{g}_{J}^{k,(1,0)}\cong\mathfrak{t}^{k+1,(1,0)}\oplus\mathfrak{g}_{J}^{k+1,(1,0)}.
\]
Then inductively,
\[
\mathfrak{g}^{(1,0)}=\mathfrak{t}^{1,(1,0)}\oplus\mathfrak{t}%
^{2,(1,0)}\oplus\cdots\oplus\mathfrak{t}^{s+1,(1,0)}.
\]
Similarly,
\[
\mathfrak{g}^{(0,1)}=\mathfrak{t}^{1,(0,1)}\oplus\mathfrak{t}%
^{2,(0,1)}\oplus\cdots\oplus\mathfrak{t}^{s+1,(0,1)}.
\]
 We remark that $\mathfrak{t}^{s+1,(1,0)}$ is indeed $\mathfrak{g}_J^{s,(1,0)}$.
\begin{proposition}\label{holomorphic}
 $\overline{\partial}\mathfrak{t}^{s+1,(1,0)}=0$. For $1\leq \ell \leq s$,
\[
\overline{\partial}\mathfrak{t}^{\ell,(1,0)}\subseteq
(\oplus_{k\leq \ell}\lie{t}^{k, *(0,1)})\wedge\lie{t}^{\ell+1, (1,0)}
\oplus \Bigl(\oplus_{k> \ell}\left( \lie{t}^{k,\ast(0,1)}\wedge\lie{t}^{k+1, (1,0)} \right)\Bigr).
\label{dbar 1-form}
\]
\end{proposition}
\bproof
Suppose that $\{\oomega_j: j=1,\dots,  \}$ is a basis for $\lieg^{*(0,1)}$ and
$\{\overline{Z}_j:j=1, \dots,  \}$ is the dual basis. For any element $V$ in
$\lieg^{1,0}$,
\[
\dbar V=\sum_j[V, {\overline{Z}}_j]^{1,0}\wedge \oomega_j.
\]
When $V$ is in $\mathfrak{t}^{s+1,(1,0)}$, it is contained in the center of the algebra
$\lieg_{\CC}$. Therefore, $\dbar V=0$.
If $V$ is in $\mathfrak{t}^{\ell,(1,0)}$ with $1\leq \ell\leq s$, by
Lemma \ref{lemma on brackets}, for all ${\overline{Z}}_j\in \liet^{k, (0,1)}$
where $k>\ell $, then
\[
[V, {\overline{Z}}_j]^{1,0}\in \liet^{k+1, (1,0)}.
\]
If $\ell\geq k$, then
$[V, {\overline{Z}}_j]^{1,0}\in \liet^{\ell+1, (1,0)}.$
\eproof

\begin{proposition}\label{simple type} When $\dim_{\CC} \liec^{1,0}=1$,
every element in $\lie{t}^{s+1,(1,0)}\wedge\lie{t}^{s,(1,0)}$
is a holomorphic Poisson structure.
\end{proposition}
\bproof  Let $C$ be an element in $\lie{t}^{s+1,(1,0)}$ and $V$ an element in
$\lie{t}^{s,(1,0)}$. Then by the previous proposition,
\[
\dbar(C\wedge V) = (\dbar C)\wedge V-C\wedge \dbar V  = -C\wedge \sum_j[V, {\overline{Z}}_j]^{1,0}
\wedge \oomega_j.
\]
By Part 2 of Lemma \ref{lemma on brackets}, for all $j$, $[V, {\overline{Z}}_j]^{1,0}$ is an element in
$\lie{t}^{s+1,(1,0)}$. When $\dim_{\CC} \liec^{1,0}=1$, every $[V, {\overline{Z}}_j]^{1,0}$
is a constant multiple of $C$. Therefore,
\[
C\wedge \sum_j[V, {\overline{Z}}_j]^{1,0}=0.
\]
It follows that $\dbar (C\wedge V)=0$. As the complex structure is abelian, Part 1 of
Lemma \ref{lemma on brackets} shows that $C\wedge V$ is Poisson.
\eproof

Consider the dual space $\mathfrak{t}^{k,\ast(0,1)}$. Since $\lie{t}^{k, (0,1)}=
\lie{g}_J^{k-1, (0,1)}/\lie{g}_J^{k, (0,1)}$, if $\oomega\in \mathfrak{t}^{k,\ast(0,1)}$,
then $\oomega ({\overline{Y}})=0$ for all $\overline{Y} \in \lie{g}^{k, (0,1)}$.

\begin{lemma}\label{semi direct}{\rm \cite{Console-Fino, Salamon}} Consider
$
\mathfrak{g}^{\ast(0,1)}=\mathfrak{t}^{1,\ast(0,1)}\oplus\mathfrak{t}
^{2,\ast(0,1)}\oplus\cdots\oplus\mathfrak{t}^{s+1,\ast(0,1)}.
$
\begin{itemize}
\item  $d\mathfrak{t}^{m,\ast(0,1)}\subseteq\left(  \oplus_{k<m}
\mathfrak{t}^{k,\ast(1,0)}\right)  \wedge\left(  \oplus_{\ell<m}
\mathfrak{t}^{\ell,\ast(0,1)}\right).$
\item $\overline{\partial}\mathfrak{t}^{k,\ast(0,1)}=0$,  for all $k$.
\item For $k\geq m$, $\left[  \mathfrak{t}^{k,(1,0)},\mathfrak{t}
^{m,\ast(0,1)}\right]  =\left\{  0\right\}   $.
\item For $k<m$, $\left[  \mathfrak{t}^{k,(1,0)},\mathfrak{t}
^{m,\ast(0,1)}\right]\subseteq\oplus_{\ell<m}\mathfrak{t}^{\ell,\ast(0,1)}$.
\end{itemize}
\end{lemma}
\bproof
Suppose that $X\in\mathfrak{t}^{k,(1,0)}$ and $\overline{Y}\in\mathfrak{t}
^{\ell,(0,1)}$ and $\overline{\omega}\in$ $\mathfrak{t}^{m,\ast(0,1)} $,
then
\[
d\overline{\omega}(X,\overline{Y})=-\overline{\omega}\left(  \left[
X,\overline{Y}\right]  \right)  .
\]
Since $\left[  \mathfrak{g}_{J}^{k,(1,0)},\mathfrak{g}_{J}^{\ell,(0,1)}\right]
\subseteq\mathfrak{g}_{J, \CC}^{1+\max{\{k,\ell\}}}$,
$\overline{\omega}\left(  \left[
X,\overline{Y}\right]  \right)=0$ except possibly when $m=1+\max{\{k,\ell\}}$.

The second item is a consequence of the first.

Since $[X, \oomega]=\iota_X d\oomega$, the third and fourth items are consequences of the first.
\eproof

\begin{corollary}\label{bracket of s n s+1}
 For all  $m$, $[  \mathfrak{t}^{s+1,(1,0)},\mathfrak{t}
^{m,\ast(0,1)}]  =0.$  For all  $m\leq s$,
$[\mathfrak{t}^{s,(1,0)},\mathfrak{t}
^{m,\ast(0,1)}]  =0.$ And
\[
[ \mathfrak{t}^{s,(1,0)},\mathfrak{t}
^{s+1,\ast(0,1)}] \subseteq
\oplus_{\ell\leq s}\mathfrak{t}^{\ell,\ast(0,1)}.
\]
\end{corollary}

\section{Computation of the map $d_2$}

\

We continue our work with the assumption that $\dim_\CC\lie{t}^{s+1,(1,0)}=1$, and
$C$ is a non-zero element in $\lie{t}^{s+1,(1,0)}$.
Let $V\in\lie{t}^{s,(1,0)}$ be non-zero, and $\Lambda=C\wedge V$.
By Proposition \ref{simple type} above,
$\Lambda$ is a holomorphic Poisson structure.
For this $\Lambda$, we now compute
$d_2: E_2^{p,q}\rightarrow E_2^{p+2, q-1}$
for all $q\geq 1$.
As a consequence of Theorem \ref{quasi isomorphic s step},
\begin{equation}
 E_2^{p,q}=
 \frac{\mbox{\rm kernel of } \adL:  H^q({\lieg}^{p,0})
 \to  H^q({\lieg}^{p+1,0})}{\mbox{\rm image of }
 \adL:  H^q({\lieg}^{p-1,0}) \to
  H^q({\lieg}^{p,0})}.
 \end{equation}


\emph{Case 1.}
Suppose that  $q=1$, $p\geq0$, and $A\in \lie{g}^{*(0,1)}\otimes \lie{g}^{p,0}$ such that
$\dbar A=0$. If $\adL(A)$ represents a zero class in
$H^1({\lieg}^{p+1,0})$, there exists $B\in \lie{g}^{p+1,0}$ such that
 \[
  {\mbox{ad}}_\Lambda(A)=\dbar B.
  \]
Since the complex structure is abelian, $\adL(B)=0$. Since $d_2([A])$ is represented by
$\adL(B)$, $d_2:E^{p,1}\to E^{p+2,0}$ is identically zero for all $p$.

\emph{Case 2.}
Consider the case when $q=2$ and $p=0$. Elements in $ \lie{g}^{*(0,2)} $ are linear combinations of $\oomega_1\wedge \oomega_2$, where $\oomega_1$ and $\oomega_2$ are elements in $\lie{t}^{m,*(0,1)}$ and $\lie{t}^{n,*(0,1)}$,
respectively.
 By the first part of Corollary \ref{bracket of s n s+1},
\[
\adL(\oomega_1\wedge \oomega_2)=C\wedge [V, \oomega_1]\wedge \oomega_2-
C\wedge [V, \oomega_2]\wedge \oomega_1.
\]
Since $V\in \lie t^{s,(1,0)}$, by Corollary
\ref{bracket of s n s+1}, $[V, \oomega_1]$ is non-zero only if $m=s+1$. However,  as
$\dim_\CC\lie{t}^{s+1,*(0,1)}=1$, by   the same corollary,
 not both $m$ and $n$ are equal to $s+1$
if $\adL(\oomega_1\wedge \oomega_2)$ is not equal to zero.
Therefore we assume that  $m=s+1$ and $n\leq s$. It follows that
\begin{equation}
\adL(\oomega_1\wedge \oomega_2)=C\wedge [V, \oomega_1]\wedge \oomega_2.
\end{equation}
Now suppose that a linear combination of such terms represents a zero class in $H^2(\lie{g}^{1,0})$, then there exists
$B$ in  $\lie{g}^{*(0,1)}\otimes \lie{g}^{1,0}$ such that
\[
\sum \adL(\oomega_1\wedge \oomega_2)=\sum C\wedge [V, \oomega_1]\wedge \oomega_2=\dbar B.
\]
Again, by Corollary \ref{bracket of s n s+1},
 $[V, \oomega_1]$ is contained in $ \oplus_{k\leq s}\lie{t}^{k, *(0,1)}$, then
\begin{equation}\label{2-forms}
\dbar B=\sum  C\wedge[V, \oomega_1]\wedge \oomega_2
\in
\oplus_n \oplus_{k\leq s}\lie{t}^{k, *(0,1)}\wedge \lie{t}^{n, *(0,1)}
\wedge \lie{t}^{s+1, (1,0)}.
 \end{equation}
 Let ${\overline\rho}\in\lie{t}^{s+1, *(0,1)}$
  be the dual of $\overline{C}$, then $B$ decomposes into the following form:
 \begin{equation}
 B=\Pi+{\overline\rho}\wedge W +\oomega \wedge C
 \end{equation}
 where $\Pi \in \left(\oplus_{k\leq s}\lie{t}^{k, *(0,1)}\right)\otimes
 \left(\oplus_{l\leq s}\lie{t}^{l, (1,0)}\right)$,
 $W\in\lie{g}^{1,0}$ and $\oomega\in \lie{g}^{*(0,1)}$.
 Then
 \[
 \dbar B=\dbar\Pi-{\overline\rho}\wedge \dbar W.
 \]
 However, from identity (\ref{2-forms}), we see that $\dbar B$ does not contain any term with
 $\overline\rho$. On the other hand, by Proposition \ref{holomorphic},
 $\dbar\Pi$ does not contribute any term in $\overline\rho$. Therefore,
 ${\overline\rho}\wedge \dbar W=0$, and
 \[
 \dbar B=\dbar \Pi.
 \]
 It follows that $d_2([A])$ is represented by $\adL(\Pi)$.
  As
 $\Pi$ is a linear combination of elements of the form
 ${\overline \rho}_k\wedge W_k$, where ${\overline\rho}_k\in
\lie{t}^{k, *(0,1)}$ for ${1\leq k\leq s}$,
and $W_k\in \oplus_{1\leq l\leq s}\lie{t}^{l,(1,0)}$, $\adL (\Pi)$ is a linear combination
$
C\wedge [V, {\overline \rho}_k]\wedge W_k.
$
However, as $V$ is in $\lie{t}^{s, (1,0)}$,
by Corollary \ref{bracket of s n s+1},
$[V, {\overline \rho}_k]\neq 0$, only when ${\overline \rho}_k\in
\lie{t}^{s+1,*(0,1)}$. Therefore, $\adL (\Pi)=0.$
So we conclude that $d_2:E^{0,2}\to E^{2,1}$ is identically zero.

\

\emph{Case 3.} Suppose that $q=2$ and $p\geq 1$. For any $A$ in $\lie{g}^{*(0,2)}\otimes \lie{g}^{(p,0)}$,
it is a linear combination of
\[
 \oomega_{m}\wedge \oomega_{n}\wedge \Theta_{m,n}
\]
where $\oomega_{m}\in\lie{t}^{m,*(0,1)}$, $\oomega_n\in \lie{t}^{n, *(0,1)}$ and
$\Theta_{m,n}\in \lie{g}^{(p,0)}$, and
$\adL (A)$ is a linear combination of
\[
C\wedge [V, \oomega_{m}\wedge \oomega_{n}]\wedge \Theta_{m,n}
=C\wedge [V, \oomega_{m}]\wedge \oomega_{n}\wedge \Theta_{m,n}
-C\wedge [V, \oomega_{n}]\wedge \oomega_{m}\wedge \Theta_{m,n}.
\]
By Corollary \ref{bracket of s n s+1}, the terms
$[V, \oomega_{m}]$ and $[V, \oomega_{n}]$ are non-zero only when $m$ and $n$ are equal to $s+1$
because $V$ is in $\liet^{s,(1,0)}$. However, as $\liet^{s+1,(1,0)}$ is only one-dimensional,
not both $\oomega_m$ and $\oomega_n$ are in $\liet^{s+1,*(0,1)}$. So $\adL(A)$ is not equal to zero only when
one of them is in $\liet^{s+1,*(0,1)}$. We assume that $\oomega_m$ spans $\liet^{s+1,*(0,1)}$. It follows that
$\adL(A)$ is a linear combination of
\[
C\wedge [V, \oomega_{s+1}]\wedge \oomega_{n}\wedge \Theta_{s+1,n},
\]
where $1\leq n\leq s$.

Now suppose that it represents a zero class in $H^2(\lie{g}^{p+1,0})$, then there exists
$B\in \lie{g}^{*(0,1)}\otimes \lie{g}^{p+1,0}$ such that
$\adL(A)=\dbar B.$
Furthermore, $B$ has the following decomposition:
\begin{equation}
 B=\Pi+{\overline\rho}\wedge W +\oomega \wedge C\wedge \Gamma
 \end{equation}
where $\Pi \in \left(\oplus_{k\leq s}\lie{t}^{k, *(0,1)}\right)\otimes
 \wedge^{p+1}\left(\oplus_{l\leq s}\lie{t}^{l, (1,0)}\right)$,
 $W\in\lie{g}^{p+1,0}$,  $\oomega\in \oplus_{k\leq s}\lie{t}^{k, *(0,1)}$, and
 $\Gamma\in \oplus_{l\leq s}\lie{t}^{l, (p,0)}$.
 Then
 \[
 \dbar B=\dbar\Pi-{\overline\rho}\wedge \dbar W+\oomega\wedge C\wedge \dbar\Gamma.
 \]
However, from Equation (\ref{2-forms}), we see that $\dbar B$ does not contain any components with
 $\overline\rho$. On the other hand, from Proposition \ref{holomorphic}, we see that
 $\dbar\Pi$ and $\dbar \Gamma$ do not contribute any terms in $\overline\rho$.
So $\adL(A)=\dbar (B)$ only if ${\overline\rho}\wedge \dbar W=0$, and
 \[
 \dbar B=\dbar \Pi+\oomega\wedge C\wedge \dbar\Gamma=\dbar (\Pi+\oomega\wedge C\wedge\Gamma).
 \]
 Then $d_2[A]$ is represented by
 \begin{equation}
 \adL(\Pi+\oomega\wedge C\wedge\Gamma)
 =C\wedge [V, \Pi+\oomega\wedge C\wedge\Gamma].
 \end{equation}

However, $V\in \lie{t}^{s, (1,0)}$,  $[V, {\overline \rho}_k]\neq 0$ only when ${\overline \rho}_k\in
\lie{t}^{s+1,*(0,1)}$. Therefore, $\adL (B)=0.$
So we conclude that for all $p\geq 0$, $d_2:E^{p,2}\to E^{p+2,1}$ is identically zero.

\

\emph{Case 4.}  Finally, we consider the case when $p\geq 0$ and $q> 2$. Let $A$
\begin{equation}
A= \oomega_{s+1}\wedge {\overline\Omega}_1\wedge \Theta_1
+{\overline\Omega}_2\wedge \Theta_2,
\end{equation}
where ${\overline\Omega}_1\in \wedge^{q-1}(\oplus_{m\leq s}\lie{t}^{m,*(0,1)})$,
${\overline\Omega}_2\in \wedge^{q}(\oplus_{m\leq s}\lie{t}^{m,*(0,1)})$,
and
$\Theta_1, \Theta_2\in \lie{g}^{p,0}$. By Corollary \ref{bracket of s n s+1},
\begin{equation}
\adL(A)=C\wedge [V, \oomega_{s+1}]\wedge{\overline\Omega}_1\wedge \Theta_1.
\end{equation}
It is contained in
\begin{equation}
\wedge^{q}(\oplus_{m\leq s}\lie{t}^{m,*(0,1)})\otimes\lie{t}^{s+1, (1,0)}\wedge \lie{g}^{p,0}.
\end{equation}
Now suppose that it represents a zero class in $H^q(\lie{g}^{p+1,0})$, then there exists
$B\in \lie{g}^{*(0,q-1)}\otimes \lie{g}^{p+1,0}$ such that
$\adL(A)=\dbar B.$
Furthermore, $B$ has the following decomposition:
\begin{equation}
 B=\Pi+{\overline\rho}\wedge {\overline\Omega}\wedge W,
 \end{equation}
where $\Pi \in \left(\wedge^{q-1}\oplus_{m\leq s}\lie{t}^{m, *(0,1)}\right)\otimes
 \lie{g}^{p+1,0}$,
 ${\overline{\Omega}}\in \wedge^{q-2}\oplus_{n\leq s}\lie{t}^{n, *(0,1)}$, and
 $W\in\lie{g}^{p+1,0}$. Then
 \begin{equation}
 \dbar B=\dbar \Pi+{\overline\rho}\wedge {\overline\Omega}\wedge \dbar W.
 \end{equation}
As $\dbar\Pi$ does not contribute any terms in ${\overline\rho}$ and such terms should be equal to zero
when $\adL(A)=\dbar (B)$, ${\overline\rho}\wedge {\overline\Omega}\wedge \dbar W=0$.
It follows that $\dbar B=\dbar \Pi$, and $d_2([A])$ is represented by $\adL(\Pi)$.
However, as
\[
\adL (\Pi)=C\wedge [V, \Pi].
\]
As a result of the semi-direct structure as seen in Lemma \ref{semi direct}, it is equal to zero.

\

Therefore, we conclude that
for all $p\geq 0, q\geq 1$, $d_2:E^{p,q}\to E^{p+2,q-1}$ is identically zero.
It concludes the proof of Theorem \ref{existence}.

\section{Examples} In \cite{CGP}, it is shown  the existence of examples of
2-step nilmanifolds    in all dimensions   {with abelian complex structures  admitting  holomorphic Poisson structures}.

In this section, we show  a sequence of high-step nilmanifolds with abelian complex structures, and
provide an explicit holomorphic Poisson structure for which the conclusion of Theorem \ref{existence}
holds. We will also provide an example of holomorphic Poisson
 structures on a complex four-dimensional  {nilmanifold} on which the holomorphic Poisson bi-complex fails
  to degenerate on its second level.

This example is inspired by the one in \cite{Milli}. The complex manifold could be considered as a tower of
 elliptic fibrations over the Kodaira surface   {\cite{BVP, Poon}}. Let $\lieg$ be a  {real}  Lie algebra with basis $\{x_1, y_1, \dots, x_n, y_n\}$ and
structure equations:
\[
[x_1, y_1]=y_2, \quad
[x_1, x_k]=[y_1, y_k]=x_{k+1}, \quad  [x_1, y_k]=-[y_1, x_k]=y_{k+1},
\]
for all $2\leq k\leq n-1$. Define an abelian complex structure by $Jx_j=y_j$ and
$Jy_j=-x_j$ for all $1\leq j\leq n$.
Let $v_j=\frac12 (x_j-iJx_j)=\frac12 (x_j-iy_j)$, then $\{v_1, \dots, v_n\}$ forms a basis for $\lieg^{1,0}$.
Let $\{\omega_1, \dots, \omega_n\}$ be the dual basis of  $\{v_1, \dots, v_n\}$,

The structure  {equations} on $\lieg_{\CC}$ with respect to the complex {basis are}
\begin{equation}\label{complexified structure eq}
[v_1, {\overline v}_1]=-\frac12(v_2-{\overline v}_2),
\quad
[v_2, {\overline v}_1]=-v_3, \quad \dots, \quad
[v_{n-1}, {\overline v}_1]=-v_n.
\end{equation}
In particular,
\begin{equation}\label{dbar v}
\dbar v_1=-\frac12 v_2\wedge \oomega_1,
\quad
\dbar v_2=-v_3\wedge\oomega_1, \quad \dots, \quad
\dbar v_{n-1}=-v_n\wedge \oomega_1,
\quad
\dbar v_n=0.
\end{equation}
It becomes apparent that $v_{n-1}\wedge v_n$ is a holomorphic bivector field.
Note that for $n\geq 5$, $\Pi_n=v_{n-3}\wedge v_n-v_{n-2}\wedge v_{n-1}$ is a non-trivial bivector field.
Moreover,
\begin{eqnarray*}
&&\dbar (v_{n-3}\wedge v_n-v_{n-2}\wedge v_{n-1})\\
&=& (\dbar v_{n-3})\wedge v_n-(\dbar v_{n-2})\wedge v_{n-1}+v_{n-2}\wedge (\dbar v_{n-1})\\
&=& -v_{n-2}\wedge\oomega_1\wedge v_n
+ v_{n-1}\wedge\oomega_1\wedge v_{n-1}-v_{n-2}\wedge  v_{n}\wedge \oomega_1=0.
\end{eqnarray*}
Thus, $\Pi_n$ determines a holomorphic Poisson structure on the nilmanifold associated to $\lieg$ when
$\dim\lieg=2n$ $(n\geq 5)$.
\

The $(1,0)$-forms $\omega_j$ satisfy the identities:
\begin{equation}
d\omega_1=0, \quad d\omega_2=\frac12 \omega_1\wedge\oomega_1, \quad
d\omega_3= \omega_2\wedge\oomega_1, \quad \cdots, \quad
d\omega_n=\omega_{n-1}\wedge\oomega_1.
\end{equation}
Taking the complex conjugation, we get
\begin{equation}\label{d oomega}
d\oomega_1=0, \quad d\oomega_2=-\frac12\omega_1\wedge\oomega_1, \quad
d\oomega_3=-\omega_1\wedge\oomega_2, \quad \cdots, \quad
d\oomega_n=-\omega_1\wedge\oomega_{n-1}.
\end{equation}
It of course follows that $\dbar\oomega_j=0$ for all $1\leq j\leq n$.
 Then we have
\begin{equation}\label{adjoint v}
[v_1, \oomega_1]=0, \quad
[v_1, \oomega_2]=-\frac12\oomega_1, \quad
[v_1, \oomega_3]=-\oomega_2, \quad
\cdots,
[v_1, \oomega_n]=-\oomega_{n-1}.
\end{equation}

It follows that when $n\geq 5$, ${\mbox \rm{ad}}_{\Pi_n}$ is identically zero. Hence the cohomology of
$\dbar_{\Pi_n}=\dbar+{\mbox\rm{ad}}_{\Pi_n}$ is equal to the cohomology of $\dbar$.

\vskip 0.5cm

Finally we discuss an example of real dimension eight, i.e. $n=4$.
We consider  $\Pi=2v_{1}\wedge v_4-v_{2}\wedge v_{3}$. Note that
\begin{eqnarray*}
&&\dbar (2v_{1}\wedge v_4-v_{2}\wedge v_{3})\\
&=& (2\dbar v_{1})\wedge v_4-(\dbar v_{2})\wedge v_{3}+v_{2}\wedge (\dbar v_{3})\\
&=& -v_{2}\wedge\oomega_1\wedge v_4
+ v_{3}\wedge\oomega_1\wedge v_{3}-v_{2}\wedge  v_{4}\wedge \oomega_1=0.
\end{eqnarray*}
Thus $\Pi=2v_{1}\wedge v_4-v_{2}\wedge v_{3}$ is indeed a holomorphic Poisson structure, and it turns out to be  more interesting.  Let us examine
${\mbox \rm{ad}}_{\Pi}$ a little further
in the next paragraphs. We take the sets of equations in (\ref{dbar v}) and (\ref{adjoint v}) for $n=4$.
The non-zero equations become
 \begin{equation}\label{dbar v when n=4}
\dbar v_1=-\frac12 v_2\wedge \oomega_1,
\quad
\dbar v_2=-v_3\wedge\oomega_1, \quad
\dbar v_{3}=-v_4\wedge \oomega_1.
\end{equation}
\begin{equation}\label{adjoint v when n=4}
[v_1, \oomega_2]=-\frac12\oomega_1, \quad
[v_1, \oomega_3]=-\oomega_2, \quad
[v_1, \oomega_4]=-\oomega_{3}.
\end{equation}
Based on the above information, we will demonstrate the following observation.
\begin{proposition} The holomorphic Poisson bi-complex associated to
$\Pi$ does not degenerate on the second level.
\end{proposition}

We will demonstrate that the map $d_2:E_2^{0,2}\to E_2^{2,1}$ is non-zero.
Recall that
\[
 E_2^{0,2}={\mbox{\rm kernel of }}\ {\mbox{\rm ad}_{\Pi}}:H^2(\CC)\to H^2(\lieg^{1,0}).
 \]
 Since all $(0,k)$-forms are $\dbar$-closed, $H^2(\CC)$ is spanned by
 $\oomega_i\wedge\oomega_j$ for all $1\leq i<j\leq 4$.
Since ${\mbox{\rm ad}_{v_2\wedge v_3}}$ is identically zero,
\begin{equation}
\ker {\mbox{\rm ad}_{\Pi}}=\ker {\mbox{\rm ad}_{v_1\wedge v_4}},
\quad
\mbox{\rm and }
\quad
{\mbox{Image}\ } {\mbox{\rm ad}_{\Pi}}={\mbox{Image}\ } {\mbox{\rm ad}_{v_1\wedge v_4}}.
\end{equation}
By (\ref{dbar v when n=4}) and (\ref{adjoint v when n=4}),
\begin{equation}
{\mbox{\rm ad}_{v_1\wedge v_4}}\oomega_2=\frac12 v_4\wedge\oomega_1
=-\frac12\dbar v_3,
\quad
{\mbox{\rm ad}_{v_1\wedge v_4}}\oomega_3=v_4\wedge \oomega_2.
\end{equation}
It follows that
\begin{eqnarray*}
&&{\mbox{\rm ad}_{\Pi}}\oomega_2\wedge \oomega_3=2{\mbox{\rm ad}_{v_1\wedge v_4}}\oomega_2\wedge \oomega_3\\
&=&2({\mbox{\rm ad}_{v_1\wedge v_4}}\oomega_2)\wedge \oomega_3
-2({\mbox{\rm ad}_{v_1\wedge v_4}}\oomega_3)\wedge \oomega_2\\
&=&\left (-\frac12\dbar v_3 \right)\wedge \oomega_3=\dbar \left(-\frac12v_3\wedge\oomega_3\right).
\end{eqnarray*}
Therefore $\oomega_2\wedge\oomega_3$ represents a non-trivial element in
$E_2^{0,2}$, and
$-2d_2(\oomega_2\wedge \oomega_3)$ is represented by
${\mbox{\rm ad}_{\Pi}}(v_3\wedge\oomega_3)$. By
(\ref{adjoint v when n=4}) it is equal to
\begin{eqnarray*}
&&2{\mbox{\rm ad}_{v_1\wedge v_4}}(v_3\wedge\oomega_3)\\
&=&-2({\mbox{\rm ad}_{v_1\wedge v_4}}\oomega_3)\wedge v_3=-2v_4\wedge[v_1, \oomega_3]\wedge v_3\\
&=&2v_4\wedge\oomega_2\wedge v_3=2v_3\wedge v_4\wedge\oomega_2.
\end{eqnarray*}
Since the image of $\dbar$ of any vector field has to have a $\oomega_1$ factor,
$v_3\wedge v_4\wedge\oomega_2$ is not $\dbar$-exact. Therefore, it represents a non-zero element in
the kernel of the map
\[
 {\mbox{\rm ad}_{\Pi}}:H^1(\lieg^{2,0})\to H^1(\lieg^{3,0}).
\]
Furthermore, for any vector $v$ in $\lieg^{1,0}$,
\[
{\mbox{\rm ad}_{v_1\wedge v_4}}(v\wedge\oomega_1)=0,
\quad
{\mbox{\rm ad}_{v_1\wedge v_4}}(v\wedge\oomega_2)=-\frac12 v\wedge v_4\wedge\oomega_1,
\]
\[
{\mbox{\rm ad}_{v_1\wedge v_4}}(v\wedge\oomega_3)=-v\wedge v_4\wedge\oomega_2,
\quad
{\mbox{\rm ad}_{v_1\wedge v_4}}(v\wedge\oomega_4)=-v\wedge v_4\wedge\oomega_3.
\]
It follows that ${\mbox{\rm ad}_{v_1\wedge v_4}}(v_3\wedge\oomega_3)=
-v_3\wedge v_4\wedge\oomega_2$. Yet $v_3\wedge\oomega_3$ is not $\dbar$-closed. Therefore,
$v_3\wedge v_4\wedge\oomega_2$ represents a non-zero element in
\[
E_2^{2,1}=
\frac{{\mbox{\rm kernel of }}\ {\mbox{\rm ad}_{\Pi}}:H^1(\lieg^{2,0})\to
H^1(\lieg^{3,0})}{{\mbox{\rm{image of }}}\ {\mbox{\rm ad}_{\Pi}}:H^1(\lieg^{1,0}))\to H^1(\lieg^{2,0})}.
\]
In other words,
$d_2(\oomega_2\wedge \oomega_3)$ is represented by a non-zero element, and hence $d_2$ is not identically
zero for the holomorphic Poisson structure
$\Pi=2v_{1}\wedge v_4-v_{2}\wedge v_{3}$.

\

\noindent{\bf Acknowledgments.} Z. Chen is partially supported by NSFC grant 11471179 and the Beijing Higher Education Young Elite
Teacher Project. A. Fino is partially supported by PRIN, FIRB  and by GNSAGA (Indam).
A. Fino and Y.S. Poon are grateful for hospitality of the Yau Mathematical Sciences Center of Tsinghua
University during their visits in summer 2014 and 2015. We would like to thank the anonymous reviewer  for very constructive and detailed comments.

\end{document}